\documentclass[showpacs,showkeys,amssymb,aps,notitlepage]{revtex4-2}

\usepackage{amsmath}

\usepackage[pdfauthor={Richard J. Mathar}, colorlinks=true,citecolor=blue,pdfpagelayout=OneColumn]{hyperref}

\usepackage{graphicx}

\DeclareMathOperator{\sgn}{sgn}

\begin{document}

\author{Richard J. Mathar}
\email{mathar@mpia-hd.mpg.de}
\homepage{https://www.mpia-hd.mpg.de/~mathar}
\affiliation{Max-Planck Institut f\"ur Astronomie, K\"onigstuhl 17, 69117 Heidelberg, Germany}

\thanks{
This work was supported by the NWO VICI grant
639.043.201
``Optical Interferometry: A new Method for Studies of Extrasolar Planets''
to A. Quirrenbach.
}

\date{\today}
\title{Geodetic Line at Constant Altitude above the Ellipsoid}

\begin{abstract}
The two-dimensional surface of a
bi-axial
ellipsoid is characterized
by the lengths of its major and minor axes.
Longitude and latitude span an angular coordinate system across.
We consider the egg-shaped surface of constant altitude
above (or below) the ellipsoid surface, and compute the
geodetic lines---lines of minimum Euclidean length---within this surface
which connect two points of fixed coordinates.
This addresses the common ``inverse'' problem of geodesics
generalized to non-zero elevations.
The system of differential equations which couples the two angular
coordinates along the trajectory is reduced to a single integral,
which is handled by Taylor expansion up to fourth
power in the eccentricity.
\end{abstract}

\maketitle

\section{Contents}
The common three parameters employed to relate Cartesian coordinates
to an ellipsoidal surface are the angles of latitude and longitude in a grid
on the surface, plus an altitude which is a shortest (perpendicular) distance to the surface.
The well-known functional relations (coordinate transformations) are summarized in
Section \ref{sec.sphero}.

The inverse problem of geodesy is to find the line embedded in the ellipsoid surface
which connects two fixed points subject to minimization of its length. We pose the
equivalent problem for lines at constant altitude, as if one would ask for the
shortest track of the center of a sphere of given radius which rolls on the ellipsoid surface
and  meets two points of the same, known altitude. In Section \ref{sec.christo},
we formulate this in terms
of the generic differential equations of geodesy parametrized by the Christoffel
symbols.

In the main Section \ref{sec.main}, the coupled system of differential equations
of latitude and longitude as a function of path length  is reduced to one degree
of freedom, here chosen to be the direction along the path at one of the fixed
terminal points, measured in the topocentric horizontal system,
and dubbed the launching angle. Closed-form expressions of this parameter
in terms of the coordinates of the terminal points have not been found; instead,
the results are presented as series expansions up to fourth power in the
eccentricity. 

The standard treatment
of this analysis is the projection
on an auxiliary sphere; this technique is (almost) completely ignored but for the 
pragmatic
aspect
that the case of zero eccentricity is a suitable zeroth-order
reference of series expansions around small eccentricities.

\section{Spheroidal Coordinates}\label{sec.sphero}
\subsection{Surface}

\begin{figure}
\includegraphics[width=8cm]{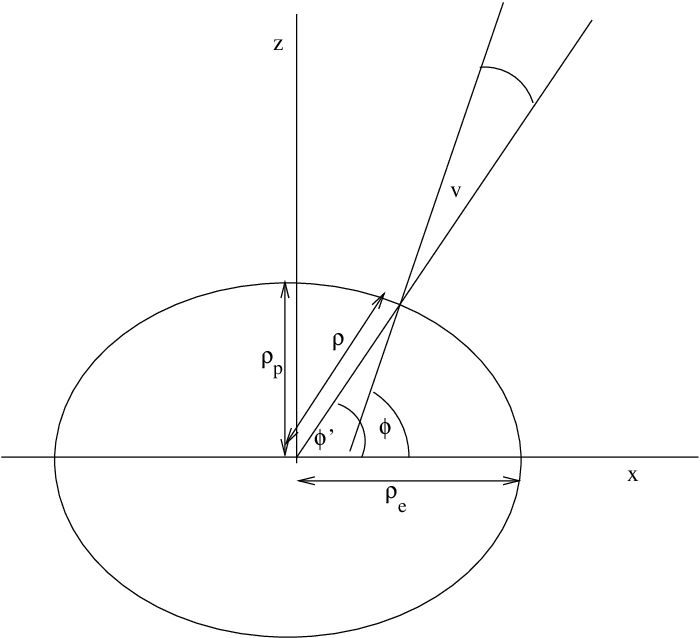}
\caption{Major semi-axis $\rho_e$, minor semi-axis $\rho_p$,
straight distance to the center of coordinates $\rho$, geocentric latitude $\phi'$,
geodetic latitude $\phi$, and their pointing difference $v$.
}
\label{fig.1}
\end{figure}

The cross section of an ellipsis of equatorial radius $\rho_e>\rho_p$
with
eccentricity $e$
in a Cartesian ($x,z$) system is:
\begin{equation}
\rho_p^2=\rho_e^2(1-e^2);
\quad
\frac{x^2}{\rho_e^2}+\frac{z^2}{\rho_p^2}=1.
\label{eq.eDef}
\end{equation}
The ellipsis defines a geocentric latitude $\phi'$ and a geodetic latitude $\phi$,
the latter measured by intersection of the normal to the tangential plane
with the equatorial plane (Figure \ref{fig.1}).
On the surface of the ellipsoid \cite[\S IX]{Smart}\cite{DufourJGeod32}\cite[\S 140]{Bouasse}:
\begin{equation}
\tan \phi=\frac{z}{x}\frac{\rho_e^2}{\rho_p^2};
\quad
\frac{z}{x}=\tan \phi'=\frac{\rho_p^2}{\rho_e^2}\tan\phi.
\end{equation}
Supposed $\rho$ denotes the distance from the center of coordinates to the point
on the geoid surface, the transformation from $(x,z)$ to $\phi$ is
\begin{equation}
x=\rho\cos\phi'=\frac{\rho_e\cos\phi}{\sqrt{1-e^2\sin^2\phi}},\quad
 z=\rho\sin\phi'=\frac{\rho_e(1-e^2)\sin\phi}{\sqrt{1-e^2\sin^2\phi}}.
\label{eq.xygeod}
\end{equation}
$v=\phi-\phi'$ defines the difference between the geodetic (astronomical) and the geocentric latitudes,
\begin{equation}
\tan v=\frac{e^2\sin(2\phi)}{2(1-e^2\sin^2\phi)}=\frac{m\sin(2\phi)}{1+m\cos(2\phi)}
=m\sin(2\phi)-m^2\sin(2\phi)\cos(2\phi)+m^3\sin(2\phi)\cos^2(2\phi)
+\cdots
;
\end{equation}
where the expansion of the denominator has been given in terms of the geometric
series of the parameter
\begin{equation}
m\equiv\frac{e^2}{2-e^2}
.
\end{equation}
The Taylor series in powers of $\sin(2\phi)$  and in powers of $m$ are
\begin{eqnarray}
v
&=&
\frac{m}{1+m}\sin(2\phi)+\frac{m^2(2+m)}{6(1+m)^3}\sin^3(2\phi)
+\frac{m^2(5+25m+15m^2+3m^3)}{40(1+m)^5}\sin^5(2\phi)+\cdots
\label{eq.vexp}
\\
&=&
\sin(2\phi)m
-\frac{1}{2}\sin(4\phi)m^2
+\frac{1}{3}\sin(6\phi)m^3
-\frac{1}{4}\sin(8\phi)m^4
+\frac{1}{5}\sin(10\phi)m^5
+\cdots
.
\end{eqnarray}
A flattening factor $f$ is also commonly defined
\cite{Kaplanarxiv06}\cite[(3.14)]{RappGG1},
\begin{equation}
f=\frac{\rho_e-\rho_p}{\rho_e}=1-\sqrt{1-e^2}
;
\quad
e^2=f(2-f).
\label{eq.flat}
\end{equation}
The reference values of the Earth ellipsoid adopted in the
WGS84 \cite{NIMA8350} are
\begin{equation}
f=1/298.257223563, \quad \rho_e=6378137.0 \text{ m}.
\end{equation}

\subsection{General Altitude}
If one moves along the direction of $\phi$ a distance $h$
away from the surface of  the geoid, the new coordinates relative to (\ref{eq.xygeod}) are
\cite{FukushimaJG79,VermeilleJG78,JonesJG76,PollardJG76,ZhangJG79,HradilekBullG50,KeelerSIR40}
\begin{equation}
x=\rho\cos\phi'+h\cos\phi;\quad
z=\rho\sin\phi'+h\sin\phi,
\end{equation}
which can be written in terms of a distance $N(\phi)$,
\begin{equation}
N(\phi)\equiv \frac{\rho_e}{\sqrt{1-e^2\sin^2\phi}}
\label{eq.Ndef}
\end{equation}
as
\begin{equation}
x=(N+h)\cos\phi;\quad
z=[N(1-e^2)+h]\sin\phi.
\label{eq.xz}
\end{equation}
\begin{figure}
\includegraphics[width=8cm]{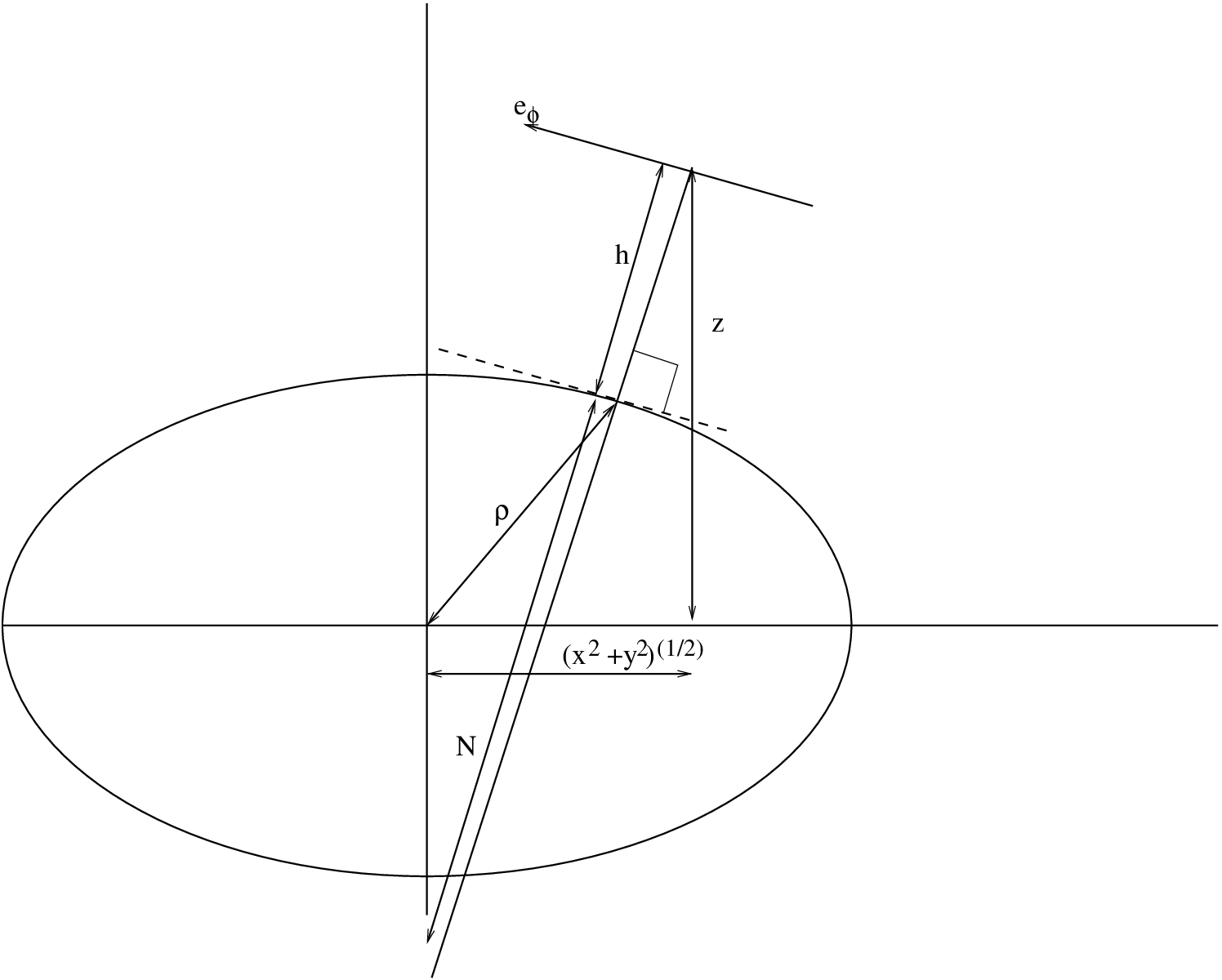}
\caption{A point with Cartesian coordinates $(x,y,z)$ has a distance
$z$ to the equatorial plane, a distance $\sqrt{x^2+y^2}$ to the polar axis,
a distance $h$ to the surface of the ellipsoid, and a distance $N+h$
to the polar axis, measured along the local normal to the surface \cite{LongCeMec12}.
The vector ${\mathbf e}_\phi$ points North at this point.
}
\end{figure}
Rotation of (\ref{eq.xz}) around the polar axis with geographic longitude
$\lambda$ defines the full 3D transformation between $(x,y,z)$ and $(\lambda,\phi,h)$,
\begin{equation}
\left(
\begin{array}{c}
x \\
y \\
z \\
\end{array}
\right)
=
\left(
\begin{array}{c}
\left[ N(\phi)+h \right] \cos \phi \cos \lambda \\
\left[ N(\phi)+h \right] \cos \phi \sin \lambda \\
\left[ N(\phi) (1-e^2)+h \right]\sin\phi \\
\end{array}
\right).
\label{eq.cart2geo}
\end{equation}

The only theme of this paper is to generalize the geodetic lines
of the literature
\cite{DufourJGeod32,SodanoJGeod32,SaitoBullGeo44,ThomasJGR70,OlanderBullGeo26,BowringBullGeo57,KarneyJG87,MaiJAG4} to the case 
of finite altitude $h\neq 0$. The physics of gravimetric or potential theory is not involved,
only the mathematics of the geometry.
It should be noted that points with constant, non-zero $h$ do \emph{not} define a surface
of an ellipsoid with effective semi-axes $\rho_{e,p}+h$---otherwise the geodetic line could
be deduced by mapping the problem onto an equivalent ellipsoidal surface \cite{DufourJGeod32}.

\section{Inverse Problem of Geodesy}\label{sec.christo}
\subsection{Topocentric Coordinate System}

A line of shortest distance at constant height $h=$const between
two points 1 and 2 is defined by minimizing the Euclidean distance,
the line integral
\begin{equation}
S=\int_1^2 \sqrt{dx^2+dy^2+dz^2}
=
\int_1^2 \sqrt{ds^2}
\equiv \int_1^2 {\cal L}ds
\end{equation}
for some parametrization $\lambda(\phi)$.
The integrand is equivalent to a Lagrange function
${\cal L}(\phi,\lambda,d\lambda/d\phi)$
or
${\cal L}(\phi,\lambda,d\phi/d\lambda)$.

The gradient of (\ref{eq.cart2geo})
with respect to  $\lambda$ and $\phi$ defines the
vectors ${\mathbf e}_{\lambda,\phi}$ that span
the topocentric tangent plane
\begin{equation}
{\bf e}_\lambda =
\left(
\begin{array}{c}
-\left[N(\phi)+h\right]\cos\phi\sin\lambda \\
\left[N(\phi)+h\right]\cos\phi\cos\lambda \\
0 \\
\end{array}
\right);
\quad
{\bf e}_\phi =
\left(
\begin{array}{c}
-\sin\phi\cos\lambda \left[M(\phi)+h\right]
   \\
-\sin\phi\sin\lambda \left[M(\phi)+h\right]
   \\
\cos\phi \left[M(\phi)+h\right]
   \\
\end{array}
\right),
\end{equation}
where a meridional radius of curvature \cite{TienstraBullGeo25}\cite[(3.87)]{RappGG1}
\begin{equation}
M(\phi)\equiv \frac{\rho_e(1-e^2)}{(1-e^2\sin^2\phi)^{3/2}}=N(\phi)\frac{1-e^2}{1-e^2\sin^2\phi}
\label{eq.Mdef}
\end{equation}
is defined to simplify the notation.
Building squares and dot products computes the three Gauss Fundamental parameters of the surface
\cite{Reichardt}
\begin{eqnarray}
{\bf e}_\lambda^2 = E &=& (N+h)^2\cos^2\phi ; \label{eq.E}\\
{\bf e}_\lambda\cdot {\bf e}_\phi = F &=& 0 ; \label{eq.F}\\
{\bf e}_\phi^2 = G \label{eq.G}
&=& 
\left[
N\frac{ 1-e^2 }{1-e^2\sin^2\phi}
+h
\right]^2
=(M+h)^2
.
\end{eqnarray}
Specializing to $h=0$ we get the You formulae \cite{YouQV}.
$E$ and $G$ provide the principal curvatures along the meridian
and azimuth \cite{DorrerQV}, and the coefficients of the metric tensor in the quadratic
form of $ds^2$,
\begin{equation}
S=\int\sqrt{Ed\lambda^2+2Fd\lambda d\phi+Gd\phi^2}
.
\label{eq.Sint}
\end{equation}

\subsection{Christoffel Symbols}

Christoffel symbols are the connection coefficients between differentials
$d{\bf e}_\epsilon$ of the topocentric axis and of the positions $d{\bf x}^\beta$, in
a generic definition
\begin{equation}
d{\bf e}_\epsilon=
\sum_{\alpha,\beta} {\bf e}_\alpha \Gamma_{\beta\epsilon}^\alpha d{\bf x}^\beta
.
\label{eq.deChr}
\end{equation}
This format is matched
by first computing the derivative of
${\bf e}_\lambda$
with respect to $\lambda$ and $\phi$ (at $h=const$):
\begin{equation}
d{\bf e}_\lambda
=
\left(
\begin{array}{c}
-(N+h)\cos\phi\cos\lambda \\ 
-(N+h)\cos\phi\sin\lambda \\ 
0 \\ 
\end{array}
\right)d\lambda
+
\left(
\begin{array}{c}
(M+h)\sin\phi\sin\lambda \\ 
-(M+h)\sin\phi\cos\lambda \\ 
0 \\ 
\end{array}
\right)d\phi
\label{eq.delam}
\end{equation}
and of
${\bf e}_\phi$
with respect to $\lambda$ and $\phi$:
\begin{equation}
d{\bf e}_\phi
=
\left(
\begin{array}{c}
(M+h)\sin\phi\sin\lambda \\ 
-(M+h)\sin\phi\cos\lambda \\ 
0 \\ 
\end{array}
\right)d\lambda
+
\left(
\begin{array}{c}
\left[Ne^2\frac{\cos^2\phi}{1-e^2\sin^2\phi}
      -3Ne^2(1-e^2)\frac{\sin^2\phi}{(1-e^2\sin^2\phi)^2}
      -(N+h)\right]\cos\phi\cos\lambda \\ 
\left[Ne^2\frac{\cos^2\phi}{1-e^2\sin^2\phi}
      -3Ne^2(1-e^2)\frac{\sin^2\phi}{(1-e^2\sin^2\phi)^2}
-(N+h)\right]\cos\phi\sin\lambda \\ 
\left[
Ne^2(1-e^2)\left(3\frac{\cos^2\phi}{(1-e^2\sin^2\phi)^2}-\frac{\sin^2\phi}{1-e^2\sin^2\phi}\right)
-[N(1-e^2)+h]
\right]\sin\phi \\ 
\end{array}
\right)d\phi
.
\label{eq.dephi}
\end{equation}
The next step splits these two equations to the expanded version of
(\ref{eq.deChr}),
\begin{eqnarray}
d{\bf e}_\lambda
&=&
({\bf e}_\lambda\Gamma_{\lambda\lambda}^\lambda +{\bf e}_\phi\Gamma_{\lambda\lambda}^\phi)d\lambda
+({\bf e}_\lambda\Gamma_{\phi\lambda}^\lambda +{\bf e}_\phi\Gamma_{\phi\lambda}^\phi)d\phi
;
\\
d{\bf e}_\phi
&=&
({\bf e}_\lambda\Gamma_{\lambda\phi}^\lambda +{\bf e}_\phi\Gamma_{\lambda\phi}^\phi)d\lambda
+({\bf e}_\lambda\Gamma_{\phi\phi}^\lambda +{\bf e}_\phi\Gamma_{\phi\phi}^\phi)d\phi
.
\end{eqnarray}
The eight $\Gamma$ are extracted by evaluating dot products of the four vector
coefficients in (\ref{eq.delam})--(\ref{eq.dephi})
by ${\bf e}_\lambda$ and ${\bf e}_\phi$,
\begin{equation}
\Gamma_{\lambda\lambda}^\lambda=
\Gamma_{\phi\lambda}^\phi
=
\Gamma_{\lambda\phi}^\phi
=
\Gamma_{\phi\phi}^\lambda
=0 ;
\label{eq.Gammazero}
\end{equation}
\begin{equation}
\Gamma_{\lambda\lambda}^\phi=(N+h)\cos\phi\sin\phi
\frac{1-e^2\sin^2\phi}{h(1-e^2\sin^2\phi)+N(1-e^2)} ;
\end{equation}
\begin{equation}
\Gamma_{\phi\lambda}^\lambda
=
\Gamma_{\lambda\phi}^\lambda
=
-\sin\phi
\frac{h(1-e^2\sin^2\phi)+N(1-e^2)}{(N+h)(1-e^2\sin^2\phi)\cos\phi} ;
\label{eq.Gammapll}
\end{equation}
\begin{equation}
\Gamma_{\phi\phi}^\phi
=
3Ne^2(1-e^2)\sin\phi\cos\phi 
\frac{1}{[h(1-e^2\sin^2\phi)+N(1-e^2)](1-e^2\sin^2\phi)}
.
\end{equation}

The Euler-Lagrange Differential Equations
$\delta\int_1^2 \sqrt{Ed\lambda^2+Gd\phi^2}=0$
for a stationary
Lagrange density $\cal L$ (at $F=0$) become the differential equations of the geodesic
\cite{Reichardt},
in the generic format
\begin{equation}
\frac{d^2 x^\epsilon}{d s^2} + \sum_{\mu\nu} \Gamma_{\mu\nu}^\epsilon
\frac{d x^\mu}{d s}\frac{d x^\nu}{d s}
=0
.
\end{equation}
The explicit write-up
\begin{equation}
\frac{d^2 \lambda}{d s^2} +
\Gamma_{\lambda\lambda}^\lambda
\frac{d \lambda}{d s}\frac{d \lambda}{d s}
+
\Gamma_{\lambda\phi}^\lambda
\frac{d \lambda}{d s}\frac{d \phi}{d s}
+
\Gamma_{\phi\lambda}^\lambda
\frac{d \phi}{d s}\frac{d \lambda}{d s}
+
\Gamma_{\phi\phi}^\lambda
\frac{d \phi}{d s}\frac{d \phi}{d s}
=0
;
\end{equation}
\begin{equation}
\frac{d^2 \phi}{d s^2} +
\Gamma_{\lambda\lambda}^\phi
\frac{d \lambda}{d s}\frac{d \lambda}{d s}
+
\Gamma_{\lambda\phi}^\phi
\frac{d \lambda}{d s}\frac{d \phi}{d s}
+
\Gamma_{\phi\lambda}^\phi
\frac{d \phi}{d s}\frac{d \lambda}{d s}
+
\Gamma_{\phi\phi}^\phi
\frac{d \phi}{d s}\frac{d \phi}{d s}
=0
,
\end{equation}
simplifies with (\ref{eq.Gammazero}) to
\begin{eqnarray}
\frac{d^2 \lambda}{d s^2} +
2\Gamma_{\lambda\phi}^\lambda
\frac{d \lambda}{d s}\frac{d \phi}{d s}
&=& 0
\label{eq.tosep}
;
\\
\frac{d^2 \phi}{d s^2} +
\Gamma_{\lambda\lambda}^\phi
\left(\frac{d \lambda}{d s}\right)^2
+
\Gamma_{\phi\phi}^\phi
\left(\frac{d \phi}{d s}\right)^2
&=& 0
\label{eq.targ}
.
\end{eqnarray}

\section{Reduction of the Differential Equations}\label{sec.main}

\subsection{Separation of Angular Variables}

Decoupling of the two differential equations (\ref{eq.tosep})--(\ref{eq.targ})
starts with the separation of variables in (\ref{eq.tosep}),
\begin{equation}
\frac{\frac{d^2\lambda}{d s^2}}{\frac{d \lambda}{d s}}
=
-2\Gamma_{\phi\lambda}^\lambda \frac{d\phi}{d s}
.
\end{equation}
Change of the integration variable on the right hand side from $s$ to $\phi$
allows to use the underivative of (\ref{eq.Gammapll})
\begin{equation}
\int \Gamma_{\lambda\phi}^\lambda (\phi)d\phi = 
\log\left\{[N(\phi)+h]\cos\phi\right\}
+const
\end{equation}
to generate a first integral
\begin{equation}
\log\frac{d \lambda}{d s}
=
-2 \log[(N+h)\cos\phi]+const
.
\end{equation}
Exponentiation yields
\begin{equation}
\frac{d \lambda}{d s}
=
c_3\frac{1}{(N+h)^2\cos^2\phi}
,
\label{eq.dlds}
\end{equation}
where
\begin{equation}
c_3\equiv 
\frac{d \lambda}{d s}_{\mid 1}
(N_1+h)^2\cos^2\phi_1
\label{eq.c3at1}
\end{equation}
is a constant for each geodesic; it plays the role of the Clairaut constant \cite{SjobergJGS2,TsengJN67},
but has length units in our case.
It has been defined with the azimuth
$\phi_1$ and 
$N_1\equiv N(\phi_1)$ at the start of the line, but could as well be
associated with any other point or the end point $2$\@.
$c_3$ is positive for trajectories starting into eastwards direction,
negative for the westwards heading, zero for routes to the poles.
Moving the square of (\ref{eq.dlds}) into (\ref{eq.targ}) yields
\begin{equation}
\frac{d^2 \phi}{d s^2}
+
\frac{\sin\phi}{(N+h)^3\cos^3\phi}
\,
\frac{1-e^2\sin^2\phi}{h(1-e^2\sin^2\phi)+N(1-e^2)}
c_3^2
+
\frac{
3Ne^2(1-e^2)\sin\phi\cos\phi 
}{[h(1-e^2\sin^2\phi)+N(1-e^2)](1-e^2\sin^2\phi)}
\left(\frac{d \phi}{d s}\right)^2
=0
.
\label{eq.dphids}
\end{equation}

To solve this differential equation, we substitute the variable
$\phi$ by its projection $\tau$ onto the polar axis,
\begin{equation}
\tau\equiv \sin\phi,
\label{eq.tauofphi}
\end{equation}
which implies transformations in the derivatives:
\begin{equation}
\frac{d\tau}{d s}=\cos\phi\frac{d\phi}{ds};
\label{eq.dtauds}
\end{equation}
\begin{equation}
\frac{d^2\tau}{d s^2}=-\sin\phi\frac{d\phi}{ds}\frac{d\phi}{ds}
+\cos\phi\frac{d^2\phi}{ds^2}
;
\end{equation}
\begin{equation}
\cos\phi\frac{d^2\phi}{ds^2}
=
\frac{d^2\tau}{d s^2}
+\tau\left(\frac{d\phi}{ds}\right)^2
=
\frac{d^2\tau}{d s^2}
+\frac{\tau}{\cos^2\phi}\left(\frac{d\tau}{ds}\right)^2
=
\frac{d^2\tau}{d s^2}
+\frac{\tau}{1-\tau^2}\left(\frac{d\tau}{ds}\right)^2
.
\label{eq.dtau2}
\end{equation}
We multiply (\ref{eq.dphids}) by $\cos\phi$, then
replace $d\phi/ds$ and $d^2\phi/ds^2$ as noted above,
\[
\frac{d^2\tau}{d s^2}
+\frac{\tau}{1-\tau^2}\left(\frac{d\tau}{ds}\right)^2
+
\frac{\tau}{(N+h)^3(1-\tau^2)}
\,
\frac{1-e^2\tau^2}{h(1-e^2\tau^2)+N(1-e^2)}
c_3^2
+
\frac{
3Ne^2(1-e^2)\tau
}{[h(1-e^2\tau^2)+N(1-e^2)](1-e^2\tau^2)}
\left(\frac{d \tau}{d s}\right)^2
=0
;
\]
\begin{multline}
\frac{d^2\tau}{d s^2}
+\frac{\tau}{1-\tau^2}
\,
\frac{
h(1-e^2\tau^2)^2
+N(1-4e^2\tau^2+2e^2+4e^4\tau^2-3e^4)
}{[h(1-e^2\tau^2)+N(1-e^2)][1-e^2\tau^2]}
\left(\frac{d\tau}{ds}\right)^2
\\
+ \frac{\tau}{(N+h)^3(1-\tau^2)}
\,
\frac{1-e^2\tau^2}{h(1-e^2\tau^2)+N(1-e^2)}
c_3^2
=0
.
\end{multline}
This is a differential equation
with no explicit appearance of the independent variable $s$,
\[
(1-\tau^2)\frac{d^2\tau}{d s^2}
+
\frac{
h(1-e^2\tau^2)^2
+N(1-e^2)(1+3e^2-4e^2\tau^2)
}{[h(1-e^2\tau^2)+N(1-e^2)][1-e^2\tau^2]}
\tau \left(\frac{d\tau}{ds}\right)^2
+ \frac{\tau c_3^2}{(N+h)^3}
\,
\frac{1-e^2\tau^2}{h(1-e^2\tau^2)+N(1-e^2)}
=0
,
\]
and the standard way of progressing is the substitution
\begin{equation}
\frac{d \tau}{ds}\equiv p;
\quad \frac{d ^2\tau}{d s^2}=p\frac{d p}{d \tau}
;
\label{eq.psubs}
\end{equation}
\begin{equation}
(1-\tau^2)p\frac{d p}{d \tau}
+
\frac{
h(1-e^2\tau^2)^2
+N(1-e^2)(1+3e^2-4e^2\tau^2)
}{[h(1-e^2\tau^2)+N(1-e^2)][1-e^2\tau^2]}
\tau p^2
+ \frac{\tau c_3^2}{(N+h)^3}
\,
\frac{1-e^2\tau^2}{h(1-e^2\tau^2)+N(1-e^2)}
=0
.
\end{equation}
This is transformed to a linear differential equation
by the further substitution $P\equiv p^2$, $d P/d\tau=2p\, dp/d\tau$,
\begin{equation}
\frac{1}{2}(1-\tau^2)\frac{d P}{d \tau}
+
\frac{
h(1-e^2\tau^2)^2
+N(1-e^2)(1+3e^2-4e^2\tau^2)
}{[h(1-e^2\tau^2)+N(1-e^2)][1-e^2\tau^2]}
\tau P
+ \frac{\tau c_3^2}{(N+h)^3}
\,
\frac{1-e^2\tau^2}{h(1-e^2\tau^2)+N(1-e^2)}
=0
.
\label{eq.PDE}
\end{equation}
The standard approach is to solve the homogeneous differential
equation first,
\begin{equation}
\frac{d P}{d \tau}
=
-
2\tau
\frac{
h(1-e^2\tau^2)^2
+N(1-e^2)(1+3e^2-4e^2\tau^2)
}{[1-\tau^2][h(1-e^2\tau^2)+N(1-e^2)][1-e^2\tau^2]}
P
.
\end{equation}
After division through $P$, the left hand side is easily integrated,
and the right hand side (incompletely) decomposed into partial fractions,
\begin{eqnarray}
\log P
&=&
- 2\int \tau
\frac{
h(1-e^2\tau^2)^2
+N(1-e^2)(1+3e^2-4e^2\tau^2)
}{[1-\tau^2][h(1-e^2\tau^2)+N(1-e^2)][1-e^2\tau^2]}
d\tau
\nonumber
\\
&=&
\int \frac{-2\tau}{1-\tau^2} d\tau
+
3\int \frac{-2e^2\tau}{1-e^2\tau^2} d\tau
+
6he^2
\int 
\frac{\tau}
{h(1-e^2\tau^2)+N(1-e^2)}
d\tau
\nonumber
\\
&=&
\log(1-\tau^2)
+3\log(1-e^2\tau^2)
-
2\log\left[h(1-e^2\tau^2)^{3/2}+\rho_e(1-e^2)\right]
+const
.
\end{eqnarray}
\[
P=
const\cdot
\frac{
(1-\tau^2)
(1-e^2\tau^2)^3
}{
\left[h(1-e^2\tau^2)^{3/2}+\rho_e(1-e^2)\right]^2
}
=
const\cdot
\frac{
(1-\tau^2)
(1-e^2\tau^2)^2
}{
\left[h(1-e^2\tau^2)+N(1-e^2)\right]^2
}
=
\frac{const\cdot (1-\tau^2)}{G(\tau)}
.
\]
Solution of the \emph{inhomogeneous} differential equation
(\ref{eq.PDE}) proceeds with the variation of the constant, the ansatz
\begin{equation}
P=
c(\tau)\cdot
\frac{
(1-\tau^2)
(1-e^2\tau^2)^2
}{
\left[h(1-e^2\tau^2)+N(1-e^2)\right]^2
}
.
\label{eq.Pofc}
\end{equation}
Back insertion
into (\ref{eq.PDE}) leads to a
first order differential equation
for $c(\tau)$,
\[
\frac{dc(\tau)}{d\tau}
\frac{1-\tau^2}{G(\tau)}
=
- \frac{2\tau c_3^2}{(1-\tau^2)(N+h)^3}
\,
\frac{1-e^2\tau^2}{h(1-e^2\tau^2)+N(1-e^2)}
,
\]
which is decomposed into partial fractions
\[
\frac{dc(\tau)}{d\tau}
=
c_3^2\left[
eN
\frac{2e\tau}{1-e^2\tau^2}
\frac{1}{(N+h)^3}
+
\frac{-2\tau}{1-\tau^2}
\frac{1}{(N+h)^2}
\right]\frac{1}{1-\tau^2}
.
\]
The ensuing integral over $d\tau$ is solved by aid of the
substitution $\tau^2=u$,
\[
c(\tau)
=
-c_3^2\frac{1}{(N+h)^2(1-\tau^2)}+const
.
\]
Back into (\ref{eq.Pofc})---using $const$
to indicate
placement of any member of an anonymous bag of constants of integration,
\begin{eqnarray}
P
&=&
\left[
const-\frac{1}{(N+h)^2(1-\tau^2)}
\right]
\frac{
c_3^2(1-\tau^2)
(1-e^2\tau^2)^2
}{
\left[h(1-e^2\tau^2)+N(1-e^2)\right]^2
}
\\
&=&
c_5
\frac{
(1-\tau^2)
(1-e^2\tau^2)^2
}{
\left[h(1-e^2\tau^2)+N(1-e^2)\right]^2
}
-
\frac{
c_3^2
(1-e^2\tau^2)^2
}{
(N+h)^2\left[h(1-e^2\tau^2)+N(1-e^2)\right]^2
}
\\
&=&
c_5
\frac{1}{(h+M)^2}
\left(
1-\tau^2
-
\frac{ c_3^2 }{ (N+h)^2 }
\right)
=c_5\frac{1-\tau^2}{G(\tau)}\left(1-\frac{c_3^2}{E(\tau)}\right)
=p^2
.
\label{eq.Pofc3}
\end{eqnarray}
The subscript $1$ denotes values at the starting point of the curve,
\begin{eqnarray}
P_1
&=&
c_5
\frac{
\cos^2\phi_1
(1-e^2 \sin^2\phi_1)^2
}{
\left[h(1-e^2\sin^2\phi_1)+N_1(1-e^2)\right]^2
}
-
\frac{
c_3^2
(1-e^2\sin^2\phi_1)^2
}{
(N_1+h)^2\left[h(1-e^2\sin^2\phi_1)+N_1(1-e^2)\right]^2
}
=p_1^2
\\
&=&
c_5
\frac{
\cos^2\phi_1
(1-e^2 \sin^2\phi_1)^2
}{
\left[h(1-e^2\sin^2\phi_1)+N_1(1-e^2)\right]^2
}
-
\frac{
(d\lambda/ds)_1^2
(N_1+h)^2\cos^4\phi_1
(1-e^2\sin^2\phi_1)^2
}{
\left[h(1-e^2\sin^2\phi_1)+N_1(1-e^2)\right]^2
}
.
\end{eqnarray}
Solving for $c_5$ yields
\begin{eqnarray}
c_5
&=&
\left(\frac{d\lambda}{ds}_{\mid 1}\right)^2
\cos^2\phi_1 (N_1+h)^2
+
\frac{
p_1^2
\left[h(1-e^2\sin^2\phi_1)+N_1(1-e^2)\right]^2
}{
\cos^2\phi_1
(1-e^2 \sin^2\phi_1)^2
}
=
\left(\frac{d\lambda}{ds}_{\mid 1}\right)^2
\cos^2\phi_1 (N_1+h)^2
+
\frac{
p_1^2
}{
\cos^2\phi_1
}
G_1
\nonumber
\\
&=&
\left(\frac{d\lambda}{ds}_{\mid 1}\right)^2
E_1
+
\left(\frac{d\phi}{ds}_{\mid 1}\right)^2
G_1
=
\frac{c_3^2}{(N_1+h)^2\cos^2\phi_1}
+
\left(\frac{d\phi}{ds}_{\mid 1}\right)^2
(M_1+h)^2
.
\label{eq.c5}
\end{eqnarray}
Compared with the differential version of (\ref{eq.Sint}),
\begin{equation}
ds^2=E d\lambda^2+G d\phi^2
;\quad
1=E \left(\frac{d\lambda}{ds}\right)^2+G \left(\frac{d\phi}{ds}\right)^2
,
\label{eq.dlsqsum}
\end{equation}
we must have $c_5=1$. Note that (\ref{eq.Pofc3}) is essentially a write-up
for $p^2\sim (d\phi/ds)^2$ and could be derived 
quickly
by inserting (\ref{eq.dlds}) directly
into (\ref{eq.dlsqsum}).

\begin{figure}
\includegraphics[width=7cm]{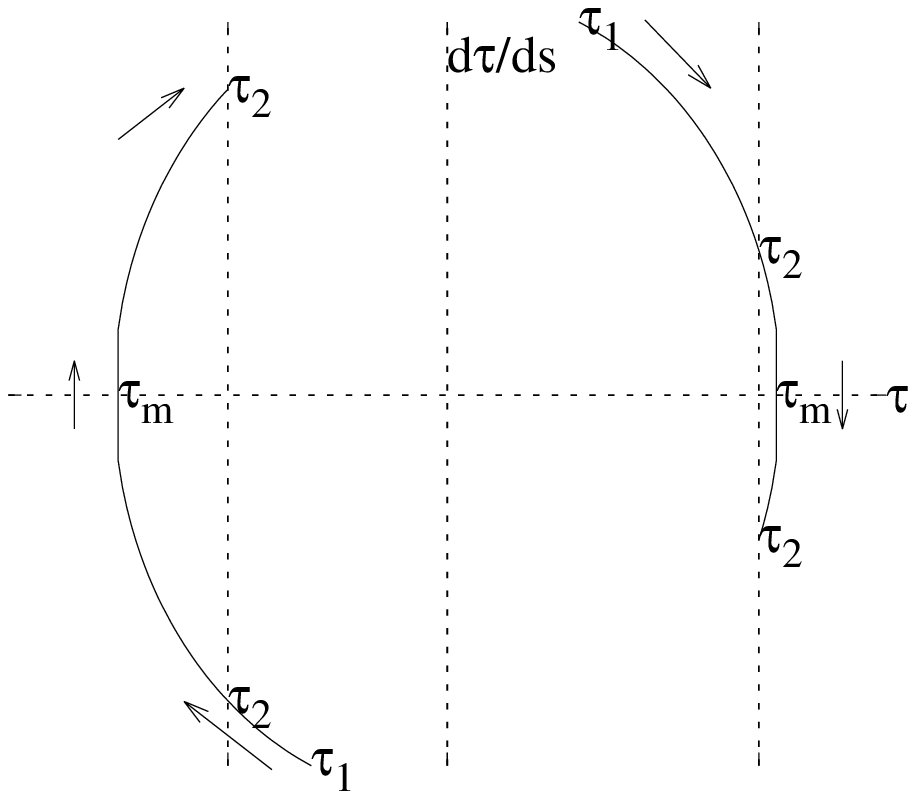}
\includegraphics[width=7cm]{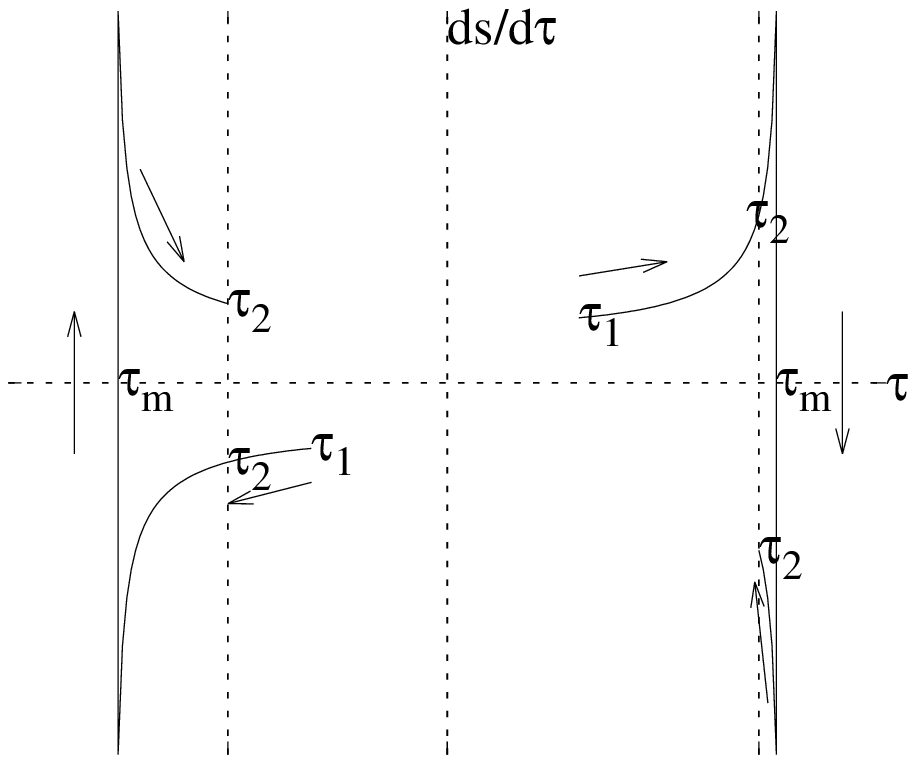}
\caption{
Two examples of trajectories of $d\tau/ds$
[left, equation (\ref{eq.dtds})] or $ds/d\tau$
[right] as a function of $\tau$. They may or may not
pass through a zero $\tau_m$ of (\ref{eq.dtds}) while connecting
the starting abscissa $\tau_1$ with the final abscissa $\tau_2$.
There are four basic topologies, depending on whether $\tau_m$
is positive or negative, and depending on whether the sign change
in $d\tau/ds$ is from $+$ to $-$ or from $-$ to $+$ at this point.
}
\label{fig.fig3}
\end{figure}

The square root of (\ref{eq.Pofc3}) is
\begin{equation}
p=\pm
\frac{
1-e^2\tau^2
}{
h(1-e^2\tau^2)+N(1-e^2)
}
\sqrt{
1-\tau^2-\frac{c_3^2}{(N+h)^2}
}
=\pm
\frac{\sqrt{1-\tau^2-\frac{c_3^2}{(N+h)^2}}}{h+M(\tau)}
=\frac{d\tau}{ds}
.
\label{eq.dtds}
\end{equation}
The sign in front of the square root
is to be chosen positive for pieces of the trajectory with 
$d\tau/ds>0$ (northern direction), negative
where $d\tau/ds<0$ (southern direction).
The value in the square root may run through zero
within one curve, $d\tau/ds=0$ at one point $\tau_m$, such that
the square root switches sign there (Figure \ref{fig.fig3}).
This happens whenever
\begin{equation}
(1-\tau_m^2)[N(\tau_m)+h]^2=c_3^2
\label{eq.dtdszero}
\end{equation}
yields a vanishing discriminant of the square
root for some $(\lambda,\phi)$, that is,
whenever the difference $\lambda_2-\lambda_1$ is sufficiently large
to create a point of minimum polar distance along the trajectory that
is not one of the terminal points.

\subsection{Launching Direction}
So far we have written the bundle of all geodesics through $(\lambda_1,\phi_1)$
in the format (\ref{eq.dtds}), which specifies the change of latitude as a function
of distance traveled. The direction at point 1 in the topocentric system of 
coordinates is represented by $c_3$.
To address the inverse problem of geodesy, that is to pick the particular
geodesics which also runs through the terminal point $(\lambda_2,\phi_2)$
in fulfillment of the Dirichlet boundary conditions,
the associated change in longitude,
some form of (\ref{eq.dlds}),
\begin{equation}
\frac{d\lambda}{ds}=c_3\frac{1}{(N+h)^2(1-\tau^2)}=\frac{c_3}{E},
\label{eq.dldsc3}
\end{equation}
must obviously get involved. 
The strategy is to write down $\tau(\lambda)$
or $\lambda(\tau)$ with $c_3$ as a parameter, then to adjust
$c_3$ to ensure
with what would be called a shooting method
that starting at point 1 at that
angle eventually passes through point 2\@.
Coupling of $\lambda$ to $\tau$ is done by division of (\ref{eq.dtds})
and (\ref{eq.dldsc3}),
\begin{equation}
\frac{d\tau}{d\lambda}=\frac{d\tau}{ds}\frac{ds}{d\lambda}
=
\frac{d\tau}{ds}/\frac{d\lambda}{ds}
=
\pm
\frac{(N+h)^2(1-\tau^2)\sqrt{1-\tau^2-\frac{c_3^2}{(N+h)^2}}}{c_3(h+M)}
.
\label{eq.dtdl}
\end{equation}
Separation of both variables yields
\begin{equation}
\pm \int_1 d\tau
\frac
{c_3(h+M)}
{(N+h)^2(1-\tau^2)\sqrt{1-\tau^2-\frac{c_3^2}{(N+h)^2}}}
=\lambda_{\mid 1}
.
\label{eq.lint}
\end{equation}

An approach of evaluating the integral by power series expansion in $h/\rho_e$
is proposed in \cite{MatharArxiv1005}.

For short paths, small $|\lambda_2-\lambda_1|$, the integral is
simply $\int_{\tau_1}^{\tau_2}$. If one of the cases occurs, where
the difference in $\lambda$ is too large for a solution, the 
additional path through the singularity $\tau_m$ is to be used [with some
local extremum in the graph $\tau(s)$], and the integral
is to be interpreted as
$\int_{\tau_1}^{\tau_m}-\int_{\tau_m}^{\tau_2}$.
Taking the sign change
and the symmetry of the integrand into account, this is
$\int_{\tau_1}^{\tau_2}+2\int_{\tau_2}^{\tau_m}=\int_{\tau_1}^{\tau_m}+\int_{\tau_2}^{\tau_m}$,
twice the underivative at $\tau_m$ minus the sum of the
underivative at $\tau_1$ and $\tau_2$.
In the right plot of Figure \ref{fig.fig3}, the difference is in including
or not including the loudspeaker shaped area between $\tau_2$ and
$\tau_m$.

$\tau_m$ is the
solution of (\ref{eq.dtdszero}), the positive value of the solution
taken if $d\tau/ds>0$ at $\phi_1$, the negative value
if $d\tau/ds<0$ at $\phi_1$.
The squared zero of $d\tau/ds$, $\tau_m^2$, is a root of the fourth-order polynomial which emerges by rewriting
(\ref{eq.dtdszero}),
\begin{multline}
e^4h^4(\tau_m^2)^4
+
2e^2h^2\left[
   \rho_e^2-h^2
   +e^2
   \left(c_3^2-h^2\right)
   \right] (\tau_m^2)^3
\\
+
\left[
\left(\rho_e^2-h^2\right)^2
+2e^2\left(
	2h^4-2c_3^2h^2-2\rho_e^2h^2-\rho_e^2c_3^2
      \right)
+e^4\left(c_3^2-h^2\right)^2
\right] (\tau_m^2)^2
\\
+
2
\left[
-\left( \rho_e^2-h^2 \right)^2
+
 c_3^2
 \left(\rho_e^2+h^2\right)
+e^2
\left(
-\left(c_3^2-h^2\right)^2
+\rho_e^2(c_3^2+h^2)
\right)
\right]
\tau_m^2
\\
+
\left[(\rho_e+h)^2-c_3^2\right]
\left[(\rho_e-h)^2-c_3^2\right]
=0.
\label{eq.taumPoly}
\end{multline}
Alternatively,
$\tau_m$
admits a Taylor expansion by
expanding the zero of
(\ref{eq.dtdszero})
in orders of $e$:
\begin{equation}
\pm \tau_m
=
 \sqrt{1-\frac{c_3^2}{H^2}}
+\frac{c_3^2 \rho_e}{2H^3}\sqrt{1-\frac{c_3^2}{H^2}}e^2
-\frac{3c_3^2 \rho_e}{8H^3}
\left[
  \left(
   \frac{\rho_e}{H}-\frac{c_3^2}{H^2}
   \right)^2-1
\right]
\sqrt{1-\frac{c_3^2}{H^2}}e^4
+O(e^6)
,
\label{eq.taumTayl}
\end{equation}
where some maximum distance
\begin{equation}
H\equiv \rho_e+h
\label{eq.Hdef}
\end{equation}
to the polar axis has been defined to condense the notation.

The integrand in (\ref{eq.lint}) has a Taylor expansion in $e$,
\begin{equation}
\pm \int d\tau
\left[
\frac{c_3/(\rho_e+h)}{(1-\tau^2)
                     \sqrt{1-\tau^2-\frac{c_3^2}{(\rho_e+h)^2}}}
-
\frac{c_3\rho_e [ (\rho_e+h)^2(2-\tau^2)-2c_3^2]
} {
 2(\rho_e+h)^4 \left(1-\tau^2-\frac{c_3^2}{(\rho_e+h)^2}\right)^{3/2}
}
e^2
+O(e^4)
\right]
=\lambda_{\mid 1}
.
\label{eq.dtdlInt}
\end{equation}
We integrate the left hand side of (\ref{eq.dtdlInt}) separately for
each power of $e$ up to $O(e^4)$,
\begin{multline}
\Big[
\arctan
 \frac{\tau\frac{c_3}{H}}
      {\sqrt{T}}
-\frac{c_3\rho_e}
 {2H^2}
\left(
 \frac{\tau}{\sqrt{T}}
+
   \arctan\frac{\tau}{\sqrt{ T }}
 \right)
e^2
\\
-\frac{c_3\rho_e}
{16 H^7}
\Big(
\frac{\tau}{T^{3/2}}
[
-2 H^2 c_3^2 \rho_e
+4 H^2 c_3^2 \rho_e T
+2 c_3^4 \rho_e
-6 H^3 c_3^2 T
+6 H^5 T
-9 H^5 T^2
-4 H^4 \rho_e T
+6 H^4 \rho_e T^2
]
\\
+
H^2
[3H^3-2\rho_e H^2-3c_3^2H+4c_3^2\rho_e]
\arctan\frac{\tau}{\sqrt{ T }}
\Big)
e^4
+O(e^6)
\Big]_{\tau_1}^{\tau_2}
=\pm (\lambda_2-\lambda_1)
,
\label{eq.c3}
\end{multline}
where
\begin{equation}
T\equiv 1-\tau^2-\left(\frac{c_3}{H}\right)^2
\label{eq.Tdef}
\end{equation}
is some convenient definition of the trajectory's distance to its
solstice $\tau_m$.
The first term is not
the principal value of the
arc tangent but its steadily defined extension through the entire
interval
of $\tau$-values, $|\tau|\le \tau_m$.
It is an odd function of $\tau$, and
phase jumps are corrected as follows: the inclination at $\tau=0$ has
(by inspection of the derivative above at $\tau=0$) the sign of $c_3$.
So whenever the triple
product $\sgn \tau \sgn c_3 \sgn (\arctan .)$ is negative, one must
shift the branch of the arctan by adding multiples of $\pi \sgn \tau \sgn c_3$
to the principal value.

Whether this is to be taken between the limits $\tau_2$
and $\tau_1$ or as the sum of two components (see above)
can be tested by integrating up to $\tau_m$, $\int_{\tau_1}^{\tau_m}$,
where the value of the underivative at $\tau_m$ is given by $(1/2)\arctan 0$,
effectively
$\frac{\pi}{2}\sgn \tau_m \sgn c_3$ after selecting the branch of the inverse
trigonometric function as described above.

Equation (\ref{eq.c3}) is solved numerically, where $c_3/H$ is the
unknown, where $\tau_{1,2}$ and $\lambda_{1,2}$ are known from the
coordinates of the two points that define the boundary value
problem, and where $e$ and $H$ are constant parameters.
The complexity of the equation suggests use of a Newton algorithm
to search for the zero, starting from (\ref{eq.c3spher}),
$c_3/H\approx \cos\phi_1\cos\phi_2 \sin(\lambda_2-\lambda_1)/\sin Z$,
as the initial estimate.

An alternative is to insert the power series
\begin{equation}
c_3=c_3^{(0)}+c_3^{(2)}e^2+c_3^{(4)}e^4+\cdots
\end{equation}
right into (\ref{eq.lint}), integrate the orders of $e$ term-by-term,
and to obtain the coefficients $c_3^{(2)}$, $c_3^{(4)}$ etc.\ by comparison
with the equivalent powers of the right hand side. $c_3^{(0)}$ is
given by (\ref{eq.c3spher});
the coefficients of the higher powers are
recursively calculated from linear equations. $c_3^{(2)}$, for example,
is determined via
\begin{equation}
\left[
\frac{\rho_e c_3^{(0)}}{H^2}
\sqrt{1-\tau^2-\left(\frac{c_3^{(0)}}{H}\right)^2}
\left[
1-\frac{{c_3^{(0)}}^2}{H^2}
\right]
\arctan\frac{H\tau}{\sqrt{1-\tau^2-\left(\frac{c_3^{(0)}}{H}\right)^2}}
+\frac{\rho_e c_3^{(0)}}{H^2}
\left[
\frac{1-{c_3^{(0)}}^2}{H^2}
\right]\tau
-2\frac{c_3^{(2)}}{H}\tau
\right]_{\tau_1}^{\tau_2}
=0
.
\end{equation}
The corresponding equation for $c_3^{(4)}$ is already too
lengthy
to be reproduced here, so no real advantage remains in comparison
with solving the non-linear equation (\ref{eq.dtdlInt}).

Once the parameter $c_3$ is known for a particular set of terminal coordinates,
$\lambda(\phi)$ is given by replacing $\tau_2$ and $\lambda_2$
in (\ref{eq.c3}) by any other generic pair of values. Two other variables
of interest along the curve, the direction and integrated distance from
the starting point, are then accessible with methods summarized in the
next, final two sub-chapters.

\subsection{Nautical course}
The nautical course at any point of the trajectory $\phi(\lambda)$ is 
the angle $\kappa$ in the topocentric coordinate system
spanned by ${\bf e}_\phi$
(direction North) and ${\bf e}_\lambda$ (direction East), measured
North over East,
\begin{equation}
\frac{d\bf r}{ds}
=
\cos\kappa \frac{{\bf e}_\phi}{|{\bf e}_\phi|}
+\sin\kappa \frac{{\bf e}_\lambda}{|{\bf e}_\lambda|}
=
\sin\kappa \frac{{\bf e}_\lambda}{\sqrt{E}}
+
\cos\kappa \frac{{\bf e}_\phi}{\sqrt{G}}
.
\end{equation}
$d{\bf r}/ds$ is the differential of (\ref{eq.cart2geo}) with respect to $s$,
where $d/ds=(d\lambda/ds)(d/d\lambda)+(d\phi/ds)(d/d\phi)$,
\begin{equation}
\frac{d\bf r}{ds}
\propto
\frac{d\lambda}{ds}{\bf e}_\lambda
+
\frac{d\phi}{ds}{\bf e}_\phi
.
\end{equation}
The sign $\propto$ indicates that the left hand side
of this equation is a vector normalized to unity, but not the right
hand side.
\begin{equation}
\frac{\sin\kappa / \sqrt{E}}{\cos \kappa/\sqrt{G}}
=
\frac{d\lambda/ds}{d\phi/ds}
=
\frac{d\lambda}{d\phi}
.
\end{equation}
Insertion of (\ref{eq.E}), (\ref{eq.G}) and (\ref{eq.dtauds}) yield
\begin{equation}
\frac{M+h}{(N+h)\cos\phi}\tan\kappa
=
\frac{1}{\frac{d\phi}{d\lambda}}
=
\frac{1}{\frac{d\phi}{d\tau}\frac{d\tau}{d\lambda}}
=
\frac{d\tau/d\phi}{\frac{d\tau}{d\lambda}}
=
\frac{\cos\phi}{\frac{d\tau}{d\lambda}}
.
\end{equation}
Mixing (\ref{eq.dtdl}) into this yields the course, supposed
$\phi$ and $c_3$ are known,
\begin{equation}
\tan\kappa
=
\pm
\,
\frac{c_3}{(N+h)\sqrt{\cos^2\phi-(\frac{c_3}{N+h})^2}}
.
\label{eq.kappa}
\end{equation}

\subsection{Distance To Terminal Points}
An implicit write-up for the distance from the origin ${\mathbf s}_1$
measured along the curve is given by separating variables in (\ref{eq.dtds}):
\begin{equation}
s=\pm
\int_{\tau_1} d\tau
\frac{
h(1-e^2\tau^2)+N(1-e^2)
}{
(1-e^2\tau^2)
\sqrt{ 1-\tau^2-\frac{c_3^2}{(N+h)^2} }
}
.
\label{eq.sIntgr}
\end{equation}

Power series expansion of the integrand in powers of $e$, then
integration term-by-term, generate a Taylor series of the form
\[
s=\pm(s^{(0)}+s^{(2)}e^2+s^{(4)}e^4+\ldots)\big|_{\tau_1}^\tau
.
\]
In the notation (\ref{eq.Tdef}), the components of the underivative read:
\begin{equation}
s^{(0)}
= H\arctan\frac{\tau}{\sqrt{1-\tau^2-c_3^2/H^2}}
= H\arctan\frac{\tau}{\sqrt{T}}
;
\end{equation}
\begin{equation}
s^{(2)}= -\frac{\rho_e}{4}
\Big[
(1+c_3^2/H^2)
\arctan\frac{\tau}{\sqrt{ T }}
+\tau
\frac{3(1-\tau^2)-c_3^2/H^2}
{\sqrt{ T }}
\Big]
;
\end{equation}
\begin{multline}
s^{(4)}=\frac{\rho_e}{64}
\Big\{
 \left[
9(c_3/H)^4
-6 (c_3/H)^2
-12 (c_3/H)^4 (\rho_e/H)
+4 (c_3/H)^2 (\rho_e/H)
-3
\right]
\arctan\frac{\tau}{\sqrt{ T }}
\\
+\frac{\tau}
{T^{3/2}}
\Big[
24 T (c_3/H)^4
-24 T (c_3/H)^2
+63 T^2 (c_3/H)^2
-8 (\rho_e/H) (c_3/H)^6
+8 (\rho_e/H) (c_3/H)^4
\\
-8 (\rho_e/H) T (c_3/H)^4
+8 (\rho_e/H) T (c_3/H)^2
-12 (\rho_e/H) T^2 (c_3/H)^2
-27 T^2
+30 T^3
\Big]
\Big\}
.
\end{multline}

\section{Summary}
Computation of the geodetic line within the iso-surface of constant altitude
above the ellipsoid
is of the same complexity as
on  its surface at zero altitude. Although
the surface is no longer an ellipsoid,
mathematics reaches an equivalent level of simplification
at which one integral is commonly
expanded in powers of the eccentricity or flattening factor.
We have done this first for the parameter which provides a solution to the
inverse problem, then for two of the basic functions, distance from
the starting point and compass course.

\appendix

\section{Reference: Spherical case}

The limit of vanishing eccentricity, $e=0$, simplifies the
curved trajectories to arcs of great circles, and presents
an easily accessible first estimate of the series expansions in orders of $e$.
The Cartesian coordinates of the two points to be connected then are
\begin{equation}
{\mathbf s}_1=(\rho_e+h)\left(
\begin{array}{c}
\cos\phi_1\cos\lambda_1\\
\cos\phi_1\sin\lambda_1\\
\sin\phi_1\\
\end{array}
\right)
;\quad
{\mathbf s}_2=(\rho_e+h)\left(
\begin{array}{c}
\cos\phi_2\cos\lambda_2\\
\cos\phi_2\sin\lambda_2\\
\sin\phi_2\\
\end{array}
\right)
.
\end{equation}
The angular separation $Z$ is derived from the dot product ${\bf s}_1\cdot {\bf s}_2 = |{\bf s}_1||{\bf s}_2|\cos Z$,
\begin{equation}
\cos Z= \sin\phi_1\sin\phi_2+\cos\phi_1\cos\phi_2\cos(\lambda_1-\lambda_2)
;\quad
0\le Z\le \pi
.
\label{eq.cosZ}
\end{equation}
Each point ${\bf s}$ on the great circle in between lies in the plane defined by the
sphere center and the terminal points, and can therefore be written as a linear
combination
\begin{equation}
{\bf s}=\alpha {\mathbf s}_1+\beta {\mathbf s}_2,\quad 0\le \alpha,\beta
.
\label{eq.sspher}
\end{equation}
The square must remain normalized to the squared radius
\begin{equation}
s^2=\alpha^2 {\mathbf s}_1^2+\beta^2 {\mathbf s}_2^2+2\alpha\beta {\mathbf s}_1\cdot {\mathbf s}_2
= \alpha^2(\rho_e+h)^2+\beta^2(\rho_e+h)^2+2\alpha\beta(\rho_e+h)^2\cos Z
,
\end{equation}
which couples the two expansion coefficients via
\begin{equation}
\alpha^2+\beta^2+2\alpha\beta\cos Z=1.
\end{equation}
One reduction to a single parameter $\xi$ to enforce
this condition is
\begin{equation}
\alpha=\frac{\cos\xi}{\sqrt{1+\sin(2\xi)\cos Z}};\quad
\beta=\frac{\sin\xi}{\sqrt{1+\sin(2\xi)\cos Z}};
\quad
0\le \xi \le \pi/2.
\end{equation}
In summary, the point (\ref{eq.sspher}) on the great circle
of radius $\rho_e+h$ between ${\mathbf s}_1$ and ${\mathbf s}_2$
has the Cartesian coordinates
\begin{eqnarray}
{\mathbf s}
&=&
(\rho_e+h)\frac{\cos\xi}{\sqrt{1+\sin(2\xi)\cos Z}}
\left(
\begin{array}{c}
\cos\phi_1\cos\lambda_1\\
\cos\phi_1\sin\lambda_1\\
\sin\phi_1\\
\end{array}
\right)
+(\rho_e+h)
\frac{\sin\xi}{\sqrt{1+\sin(2\xi)\cos Z}}
\left(
\begin{array}{c}
\cos\phi_2\cos\lambda_2\\
\cos\phi_2\sin\lambda_2\\
\sin\phi_2\\
\end{array}
\right)
\label{eq.xi}
\\
& \equiv &
(\rho_e+h)
\left(
\begin{array}{c}
\cos\phi(\xi)\cos\lambda(\xi)\\
\cos\phi(\xi)\sin\lambda(\xi)\\
\sin\phi(\xi)\\
\end{array}
\right)
.
\end{eqnarray}
The $z$-component of this,
\begin{equation}
\sin\phi(\xi)=\frac{\cos\xi\sin\phi_1+\sin\xi\sin\phi_2}{\sqrt{1+\sin(2\xi)\cos Z}}
,
\label{eq.phixi}
\end{equation}
allows to convert $\xi$ into $\varphi$.
No ambiguity with respect to the branch
of the $\arcsin$ arises since $-\pi/2\le\phi\le\pi/2$.
The ratio of the $y$ and $x$-components demonstrates the dependence of $\lambda$ on $\xi$,
\begin{equation}
\tan\lambda(\xi)=
\frac{\cos\xi\cos\phi_1\sin\lambda_1+\sin\xi\cos\phi_2\sin\lambda_2}
{\cos\xi\cos\phi_1\cos\lambda_1+\sin\xi\cos\phi_2\cos\lambda_2}
,
\label{eq.lambdaxi}
\end{equation}
where the $\arctan$ branch is defined by considering separately the numerator
and denominator under the restrictions that no signs are canceled.
The azimuth angle $\sigma$ between the points at $\xi=0$ and at $\xi$
follows from
${\mathbf s}\cdot {\mathbf s}_1 = (\rho_e+h)^2\cos\sigma$,
\begin{eqnarray}
\cos\sigma 
&=& 
\left(
\frac{\cos\xi}{\sqrt{1+\sin(2\xi)\cos Z}}
\left(
\begin{array}{c}
\cos\phi_1\cos\lambda_1\\
\cos\phi_1\sin\lambda_1\\
\sin\phi_1\\
\end{array}
\right)
+
\frac{\sin\xi}{\sqrt{1+\sin(2\xi)\cos Z}}
\left(
\begin{array}{c}
\cos\phi_2\cos\lambda_2\\
\cos\phi_2\sin\lambda_2\\
\sin\phi_2\\
\end{array}
\right)
\right)
\cdot
\left(
\begin{array}{c}
\cos\phi_1\cos\lambda_1\\
\cos\phi_1\sin\lambda_1\\
\sin\phi_1\\
\end{array}
\right)
\nonumber
\\
&=&  
\frac{\sin \xi\cos Z+\cos\xi}
{\sqrt{1+\sin(2 \xi)\cos Z}}
;\quad
0\le \sigma \le Z
.
\label{eq.cossig}
\end{eqnarray}
The path length along the great circle perimeter is simply the
radial distance $\rho_e+h$ to the center of coordinates times the azimuth angle
$\sigma$ measured in radians,
\begin{equation}
s=\int_1 \sqrt{ds^2}=(\rho_e+h)\sigma=(\rho_e+h)\arccos
\frac{\sin \xi\cos Z+\cos\xi}
{\sqrt{1+\sin(2 \xi)\cos Z}}
;\quad
0\le s\le (\rho_e+h) Z
.
\label{eq.sigma}
\end{equation}
To calculate the direction in the
${\mathbf e}_\lambda$-
${\mathbf e}_\phi$-plane
at the starting point ${\mathbf s}_1$, we employ $\xi$ as the
parameter that  mediates between $\lambda$ and $\phi$:
\begin{equation}
\frac{d\lambda}{ds}
=
\frac{1}{\rho_e+h}
\frac{d\lambda}{d\sigma}
=
\frac{1}{\rho_e+h}
\frac{d\lambda}{d\xi}\,\frac{d\xi}{d\sigma}
=
\frac{1}{\rho_e+h}
\frac{d\lambda}{d\xi}
/
\frac{d\sigma}{d\xi}
.
\label{eq.dlambds}
\end{equation}
To calculate $d\lambda/d\xi$, use the derivative of
(\ref{eq.lambdaxi}) with respect to $\xi$,
\begin{eqnarray}
\frac{1}{\cos^2\lambda}\frac{d\lambda}{d\xi}
&=&
\frac{-\sin\xi\cos\phi_1\sin\lambda_1+\cos\xi\cos\phi_2\sin\lambda_2}
{\cos\xi\cos\phi_1\cos\lambda_1+\sin\xi\cos\phi_2\cos\lambda_2}
\nonumber
\\
&&
-
(\cos\xi\cos\phi_1\sin\lambda_1+\sin\xi\cos\phi_2\sin\lambda_2)
\frac{-\sin\xi\cos\phi_1\cos\lambda_1+\cos\xi\cos\phi_2\cos\lambda_2}
{(\cos\xi\cos\phi_1\cos\lambda_1+\sin\xi\cos\phi_2\cos\lambda_2)^2}
.
\end{eqnarray}
In particular at the starting point, where $\xi=0$,
$\lambda=\lambda_1$ and
$\phi=\phi_1$,
\[
\frac{1}{\cos^2\lambda_1}\frac{d\lambda}{d\xi}_{\mid 1}
=
\frac{\cos\phi_2\sin\lambda_2}
{\cos\phi_1\cos\lambda_1}
-
\cos\phi_1\sin\lambda_1
\frac{\cos\phi_2\cos\lambda_2}
{(\cos\phi_1\cos\lambda_1)^2}
.
\]
By multiplication with $\cos\phi_1 \cos^2\lambda_1$
\begin{equation}
\cos\phi_1\frac{d\lambda}{d\xi}_{\mid 1}
=
\cos\phi_2\sin\lambda_2
\cos\lambda_1
-
\sin\lambda_1
\cos\phi_2\cos\lambda_2
=
\cos\phi_2
\sin(\lambda_2-\lambda_1)
.
\label{eq.dlambdxi}
\end{equation}
To calculate $d\sigma/d\xi$, we convert the cosine in (\ref{eq.cossig}) to the sine,
\begin{eqnarray}
\sin\sigma
&=&
\sqrt{1-\cos^2\sigma}
=
\frac{ \sin\xi\sin Z }{ \sqrt{1+\sin(2\xi)\cos Z} }
.
\end{eqnarray}
The derivative of this with respect to $\xi$ is
\[
\cos\sigma \frac{d\sigma}{d\xi}
=
\cos\xi\sin Z
\frac{1}{\sqrt{1+\sin(2\xi)\cos Z}}
-\frac{1}{2}
\sin\xi\sin Z\frac{2\cos(2\xi)\cos Z}{(1+\sin(2\xi)\cos Z)^{3/2}}
,
\]
in particular at the starting point, $\sigma=\xi=0$,
\begin{equation}
\frac{d\sigma}{d\xi}_{\mid 1}
=
\sin Z
.
\end{equation}
Insert this derivative and (\ref{eq.dlambdxi})
back into (\ref{eq.dlambds})
\begin{equation}
\frac{d\lambda}{ds}_{\mid 1}
=
\frac{1}{\rho_e+h}
\frac{\cos\phi_2\sin(\lambda_2-\lambda_1)}{\cos\phi_1}
\frac{1}{\sin Z}
,
\label{eq.dlds1}
\end{equation}
to obtain the master parameter $c_3$ for the spherical case
with (\ref{eq.c3at1}),
\begin{equation}
c_3^{(0)}
=
(\rho_e+h)
\frac{\cos\phi_1\cos\phi_2\sin(\lambda_2-\lambda_1)}{\sin Z}
;\quad (e=0)
.
\label{eq.c3spher}
\end{equation}

\section{Notations } 
\begin{tabular}{ll}
$c_5, c_3$        & constants of integration\\
$E, E_1$        & Gauss parameter eq.\ (\ref{eq.E}) and its value at curve origin $(\lambda_1,\phi_1,h)$\\
$E(.\setminus.)$   & Incomplete Elliptic Integral of the Second Kind in the Abramowitz-Stegun notation \cite[\S 17]{AS}\\
$e$        & eccentricity of the ellipsoid, eq.\ (\ref{eq.eDef})\\
${\mathbf e}_{\lambda,\phi}$        & horizontal topocentric coordinate vectors at ($\lambda,\phi,h$)\\
$f$        & flattening, eq.\ (\ref{eq.flat})\\
$F$        & Gauss parameter, eq.\ (\ref{eq.F}) \\
$G, G_1$        & Gauss parameter, eq.\ (\ref{eq.G}) and its value at curve origin \\
$\Gamma_{..}^.$   & Christoffel symbols in $(\lambda,\phi,h)$ coordinates\\
$h$        & vertical distance to surface of ellipsoid \\
$H$        & maximum distance to polar axis (in the equatorial plane), eq.\ (\ref{eq.Hdef})\\
$\kappa$        & nautical angle, North over East in the topocentric tangential plane \\
$\lambda$        & longitude\\
$M$        & a radius of curvature on the ellipsoid surface, eq.\ (\ref{eq.Mdef})\\
$N$        & a radius of curvature on the ellipsoid surface, eq.\ (\ref{eq.Ndef}) \\
$\rho$        & distance ellipsoid center to foot point on the surface\\
$\rho_{e,p}$        & equatorial, polar radius of ellipsoid, eq.\ (\ref{eq.eDef})\\
$s$        & distance along geodetic line measured from curve origin \\
${\mathbf s}$        & vector from ellipsoid center to point on geodetic line\\
$S$        & length of curved geodetic trajectory; equals $s_2$ at $(\lambda_2,\phi_2,h)$\\
$\sgn$        & sign function, $\pm 1$ or $0$ \\
$\sigma$        & azimuth along great circle, eq.\ (\ref{eq.sigma})\\
$T$        & normalized distance from closest polar approach, eq (\ref{eq.Tdef}) \\
$\tau, \tau_1$        & $\sin \phi$ and its value at start of curve, eq.\ (\ref{eq.tauofphi})\\
$v$        & geodetic minus geocentric latitude\\
$\phi$        & geodetic latitude\\
$\phi'$        & geocentric latitude\\
$x,y,z$        & Cartesian coordinates from ellipsoid center, eq.\ (\ref{eq.cart2geo})\\
$\xi$        & parametrization of great circle (spherical case), eq.\ (\ref{eq.xi})\\
$Z$        & cone angle of circular section (spherical case), eq.\ (\ref{eq.cosZ})\\
\end{tabular}

\section{FEM implementation}
\subsection{Compilation}

The simplest numerical solution of the inverse problem
without restriction on
the eccentricity could
be a finite-element (FEM) integration
of (\ref{eq.lint}) and iterative adjustment of the single
free parameter $c_3$. The Java program in the \texttt{anc} directory implements this approach.
It is either compiled with
\begin{verbatim}
cd anc ; make
\end{verbatim}
on Linux systems or by manually executing the compiler steps in \texttt{anc/Makefile} .
It is then called  with
\begin{itemize}
\item
\begin{verbatim}
java -jar Geod.jar 
\end{verbatim}
which will query the parameters interactively 
\item
or providing all parameters with

\texttt{java -jar Geod.jar [-e } $e$
\texttt{] [-R }$\rho_e$
\texttt{] [-h }\textit{h}
\texttt{] [-s } $N_s$
\texttt{] [-u }$N_2$\texttt{]} $\phi_1$ $\lambda_1$ $\phi_2$ $\lambda_2$

on the command line.
The last four parameters are the geodetic angles of the start and final point in degrees.
$N_s$ is the positive integer number of finite elements in the interval $\lambda$-interval
in the approximation, and $N_2$ in the range 1\ldots $N_s$  (and typically a divisor of $N_s$)
is the small positive integer subsampling number of the results
along the trajectory which are actually printed to the standard output.
The square brackets above indicate that the switches are optional because they have 
default values; the square brackets are not part of the command line syntax.
\item
or calling the program as

\texttt{java -cp Geod.jar org.nevec.rjm.GeodGUI}

which allows to type the parameters in fields of a graphical user interface (GUI);  
the purpose of that wrapper is to embed the program in web applications
that support the Java Network Launch Protocol (JNLP).
\end{itemize}

The output lines contain the information on the trajectory: the first three values
are the $x$, $y$ and $z$ Cartesian coordinates of the subsampled points, the next three
values are $\phi$, $\lambda$ and $\kappa$ at these points in degrees, and the last value is
the length $s$ travelled from the start point.

\subsection{Overview of Functions}

Member functions in overview: the constructors \texttt{Geod} define a surface
from the parameters $\rho_e$, $e$ and $h$
in which the geodetic line is embedded.
\texttt{getCartesian} computes the vector (\ref{eq.cart2geo}).
\texttt{curvN} computes (\ref{eq.Ndef}).
\texttt{dNdtau} computes
\begin{equation}
\frac{dN}{d\tau}=N(\tau)\frac{e^2\tau}{1-e^2\tau^2}.
\end{equation}
\texttt{flatt} computes (\ref{eq.flat}).
\texttt{curvM} computes (\ref{eq.Mdef}).
\texttt{dMdtau} computes
\begin{equation}
\frac{dM}{d\tau}= 3M(\tau)\frac{e^2\tau}{1-e^2\tau^2}.
\end{equation}
\texttt{d2Mdtau2} computes
\begin{equation}
\frac{d^2M}{d\tau^2}= 3M(\tau)\frac{e^2(1+4e^2\tau^2)}{(1-e^2\tau^2)^2}.
\end{equation}
\texttt{GaussE} computes (\ref{eq.E}).
\texttt{dEdtau} computes
\begin{equation}
\frac{dE}{d\tau}=-2\tau(N+h)(M+h).
\end{equation}
\texttt{d2Edtau2} computes the next higher derivative,
$ d^2E/d\tau^2 $.
\texttt{GaussG} computes (\ref{eq.G}).
\texttt{dtaudlambda} computes $d\tau/d\lambda$ via (\ref{eq.dtdl}),
referencing one factor to (\ref{eq.dtds}).
\texttt{dsdlambda} is $ds/d\lambda$, the inverse of (\ref{eq.dldsc3}).
\texttt{discrT} computes $1-\tau^2-c_3^2/(N+h)^2$, which generalizes
(\ref{eq.Tdef}) to nonzero $e$. Its derivatives
\begin{eqnarray}
\frac{d}{d\tau}\left[1-\tau^2-\frac{c_3^2}{[N(\tau)+h]^2}\right]
&=&
2\tau\left[\frac{Ne^2c_3^2}{(N+h)^3(1-e^2\tau^2)}-1\right];\\
\frac{d^2}{d\tau^2}\left[1-\tau^2-\frac{c_3^2}{[N(\tau)+h]^2}\right]
&=&
-2\left[1-\frac{Nc_3^2e^2}{(N+h)^3(1-e^2\tau^2)}
       \left(1+\frac{3he^2\tau^2}{(N+h)(1-e^2\tau^2)}\right)\right],
\end{eqnarray}
are implemented in \texttt{dTdtau} and \texttt{d2Tdtau2}.
\texttt{dtauds} calculates (\ref{eq.dtds}).
\texttt{d2taudlambda2} calculates the derivative of (\ref{eq.dtdl}),
\begin{equation}
\frac{d^2\tau}{d\lambda^2}
=
\frac{d}{d\lambda}\frac{E(\tau)\sqrt{1-\tau^2-c_3^2/(N+h)^2}}{c_3(h+M)}
=
\frac{d\tau}{d\lambda}\frac{d}{d\tau}\frac{E(\tau)\sqrt{1-\tau^2-c_3^2/(N+h)^2}}{c_3(h+M)}
.
\end{equation}
\texttt{d3taudlambda3} is the next higher order application
of Bruno di Fa\`a's formula to relegate derivatives $d/d\lambda$
to derivatives $d/d\tau$,
\cite[0.430]{GR},
\begin{equation}
\frac{d^3\tau}{d\lambda^3}
=
\frac{d^2\tau}{d\lambda^2}
\frac{d}{d\tau}\frac{E(\tau)\sqrt{1-\tau^2-c_3^2/(N+h)^2}}{c_3(h+M)}
+
\left( \frac{d\tau}{d\lambda}\right)^2
\frac{d^2}{d\tau^2}\frac{E(\tau)\sqrt{1-\tau^2-c_3^2/(N+h)^2}}{c_3(h+M)}
.
\end{equation}
\texttt{d2sdlambda2} calculates
\begin{equation}
\frac{d^2s}{d\lambda^2}=\frac{d}{d\lambda}\frac{E}{c_3}
=
\frac{1}{c_3}\,\frac{d\tau}{d\lambda}\,\frac{d}{d\tau}E
.
\end{equation}
\texttt{d3sdlambda3} calculates
\begin{equation}
\frac{d^3s}{d\lambda^3}=\frac{d^2}{d\lambda^2}\frac{E}{c_3}
=
\frac{1}{c_3}\left[(\frac{d\tau}{d\lambda})^2\,\frac{d^2}{d\tau^2}E
+\frac{d^2\tau}{d\lambda^2}\,\frac{d}{d\tau}E
\right]
.
\end{equation}
\texttt{c3Sphere} returns the estimate (\ref{eq.c3spher}).
\texttt{dtaudsSignum} returns the sign of $d\tau/ds$ at $\tau_1$, obtained
by considering the sign of the derivative of
(\ref{eq.phixi}) with respect to $\xi$.
\texttt{adjLambdaEnd} modifies $\phi_2$ modulo $2\pi$ to select the
smallest value of $|\phi_2-\phi_1|$.
\texttt{nautAngle} computes $\kappa$ from (\ref{eq.kappa}).

\texttt{tauShoot} walks along a geodetic line on a discrete mesh of
width $\Delta\lambda$ by extrapolating
\begin{equation}
\tau_{\lambda+\Delta\lambda}\approx \tau_{\lambda}+\frac{d\tau}{d\lambda}\Delta\lambda
+\frac{1}{2}\frac{d^2\tau}{d\lambda^2}(\Delta\lambda)^2
+\frac{1}{6}\frac{d^3\tau}{d\lambda^3}(\Delta\lambda)^3
+\frac{1}{24}\frac{d^4\tau}{d\lambda^4}(\Delta\lambda)^4,
\end{equation}
initialized at $\lambda_1,\phi_1$, given $c_3$. The equivalent formula is
used to
build up $s_{\lambda+\Delta\lambda}$. 
\texttt{c3shoot} calls \texttt{tauShoot} four times to adjust
$c_3$ such that the error by which $\phi_2$ was missed---returned by \texttt{tauShoot}---is minimized. The first call assumes (\ref{eq.c3spher}),
the second takes an arbitrary small offset, and the third and fourth
estimates are from linear and quadratic interpolations in the earlier
calls to zoom into a root of this error as a function of $c_3$.
The last of these runs tabulates the Cartesian coordinates (\ref{eq.cart2geo}),
$\lambda$, $\phi$, $\kappa$ and $s$ on a subgrid of the $\lambda$-mesh.
\texttt{main} collects some adjustable parameters plus the pairs
$(\lambda_1,\phi_1)$ and $(\lambda_2,\phi_2)$, and calls \texttt{c3shoot}
to solve the inverse problem of geodetics.

\bibliographystyle{apsrmp}

\bibliography{all}

\begin{thebibliography}{30}
\expandafter\ifx\csname natexlab\endcsname\relax\def\natexlab#1{#1}\fi
\expandafter\ifx\csname bibnamefont\endcsname\relax
  \def\bibnamefont#1{#1}\fi
\expandafter\ifx\csname bibfnamefont\endcsname\relax
  \def\bibfnamefont#1{#1}\fi
\expandafter\ifx\csname citenamefont\endcsname\relax
  \def\citenamefont#1{#1}\fi
\expandafter\ifx\csname url\endcsname\relax
  \def\url#1{\texttt{#1}}\fi
\expandafter\ifx\csname urlprefix\endcsname\relax\def\urlprefix{URL }\fi
\providecommand{\bibinfo}[2]{#2}
\providecommand{\eprint}[2][]{\url{#2}}

\bibitem[{\citenamefont{Abramowitz and Stegun}(1972)}]{AS}
\bibinfo{editor}{\bibnamefont{Abramowitz}, \bibfnamefont{M.}}, and
  \bibinfo{editor}{\bibfnamefont{I.~A.} \bibnamefont{Stegun}} (eds.),
  \bibinfo{year}{1972}, \emph{\bibinfo{title}{Handbook of Mathematical
  Functions}} (\bibinfo{publisher}{Dover Publications}, \bibinfo{address}{New
  York}), \bibinfo{edition}{9th} edition, ISBN \bibinfo{isbn}{0-486-61272-4}.

\bibitem[{\citenamefont{Bouasse}(1919)}]{Bouasse}
\bibinfo{author}{\bibnamefont{Bouasse}, \bibfnamefont{H.}},
  \bibinfo{year}{1919}, \emph{\bibinfo{title}{Geographie Math\'ematique}}
  (\bibinfo{publisher}{Librairie Delagrave}, \bibinfo{address}{Paris}).

\bibitem[{\citenamefont{Bowring}(1983)}]{BowringBullGeo57}
\bibinfo{author}{\bibnamefont{Bowring}, \bibfnamefont{B.~R.}},
  \bibinfo{year}{1983}, \bibinfo{journal}{Bull. Geod.}
  \textbf{\bibinfo{volume}{57}}(\bibinfo{number}{1--4}), \bibinfo{pages}{109}.

\bibitem[{\citenamefont{Dorrer}(1999)}]{DorrerQV}
\bibinfo{author}{\bibnamefont{Dorrer}, \bibfnamefont{E.}},
  \bibinfo{year}{1999}, in \emph{\bibinfo{booktitle}{Quo vadis
  geodesia\ldots?}}, edited by
  \bibinfo{editor}{\bibfnamefont{F.}~\bibnamefont{Krumm}} and
  \bibinfo{editor}{\bibfnamefont{V.~S.} \bibnamefont{Schwarze}}
  (\bibinfo{publisher}{Universit\"at Stuttgart}), Technical Reports of Geodesy
  and GeoInformatics, ISSN \bibinfo{issn}{0933-2839}.

\bibitem[{\citenamefont{Dufour}(1958)}]{DufourJGeod32}
\bibinfo{author}{\bibnamefont{Dufour}, \bibfnamefont{H.~M.}},
  \bibinfo{year}{1958}, \bibinfo{journal}{Bull. Geod.}
  \textbf{\bibinfo{volume}{32}}(\bibinfo{number}{2}), \bibinfo{pages}{26}.

\bibitem[{\citenamefont{Fukushima}(2006)}]{FukushimaJG79}
\bibinfo{author}{\bibnamefont{Fukushima}, \bibfnamefont{T.}},
  \bibinfo{year}{2006}, \bibinfo{journal}{J. Geod.}
  \textbf{\bibinfo{volume}{79}}(\bibinfo{number}{12}), \bibinfo{pages}{689}.

\bibitem[{\citenamefont{Gradstein and Ryshik}(1981)}]{GR}
\bibinfo{author}{\bibnamefont{Gradstein}, \bibfnamefont{I.}}, and
  \bibinfo{author}{\bibfnamefont{I.}~\bibnamefont{Ryshik}},
  \bibinfo{year}{1981}, \emph{\bibinfo{title}{Summen-, {P}rodukt- und
  {I}ntegraltafeln}} (\bibinfo{publisher}{Harri Deutsch},
  \bibinfo{address}{Thun}), \bibinfo{edition}{1st} edition, ISBN
  \bibinfo{isbn}{3-87144-350-6}.

\bibitem[{\citenamefont{Hradilek}(1976)}]{HradilekBullG50}
\bibinfo{author}{\bibnamefont{Hradilek}, \bibfnamefont{L.}},
  \bibinfo{year}{1976}, \bibinfo{journal}{Bull. Geod.}
  \textbf{\bibinfo{volume}{50}}(\bibinfo{number}{4}), \bibinfo{pages}{301}.

\bibitem[{\citenamefont{Jones}(2004)}]{JonesJG76}
\bibinfo{author}{\bibnamefont{Jones}, \bibfnamefont{G.~C.}},
  \bibinfo{year}{2004}, \bibinfo{journal}{J. Geod.}
  \textbf{\bibinfo{volume}{76}}(\bibinfo{number}{8}), \bibinfo{pages}{437}.

\bibitem[{\citenamefont{Kaplan}(2006)}]{Kaplanarxiv06}
\bibinfo{author}{\bibnamefont{Kaplan}, \bibfnamefont{G.~H.}},
  \bibinfo{year}{2006}, \bibinfo{journal}{arXiv:astro-ph/0602086}
  \eprint{astro-ph/0602086}.

\bibitem[{\citenamefont{Karney}(2013)}]{KarneyJG87}
\bibinfo{author}{\bibnamefont{Karney}, \bibfnamefont{C.~F.~F.}},
  \bibinfo{year}{2013}, \bibinfo{journal}{J. Geod.}
  \textbf{\bibinfo{volume}{87}}, \bibinfo{pages}{43}.

\bibitem[{\citenamefont{Keeler and Nievergelt}(1998)}]{KeelerSIR40}
\bibinfo{author}{\bibnamefont{Keeler}, \bibfnamefont{S.~P.}}, and
  \bibinfo{author}{\bibfnamefont{Y.}~\bibnamefont{Nievergelt}},
  \bibinfo{year}{1998}, \bibinfo{journal}{SIAM Review}
  \textbf{\bibinfo{volume}{40}}(\bibinfo{number}{2}), \bibinfo{pages}{300}.

\bibitem[{\citenamefont{Long}(1975)}]{LongCeMec12}
\bibinfo{author}{\bibnamefont{Long}, \bibfnamefont{S.~A.~T.}},
  \bibinfo{year}{1975}, \bibinfo{journal}{Cel.\ Mech.}
  \textbf{\bibinfo{volume}{12}}(\bibinfo{number}{2}), \bibinfo{pages}{225}.

\bibitem[{\citenamefont{Mai}(2010)}]{MaiJAG4}
\bibinfo{author}{\bibnamefont{Mai}, \bibfnamefont{E.~.}}, \bibinfo{year}{2010},
  \bibinfo{journal}{J. Appl. Geodesy} \textbf{\bibinfo{volume}{4}},
  \bibinfo{pages}{145}.

\bibitem[{\citenamefont{Mathar}(2010)}]{MatharArxiv1005}
\bibinfo{author}{\bibnamefont{Mathar}, \bibfnamefont{R.~J.}},
  \bibinfo{year}{2010}, \bibinfo{journal}{arXiv:1005.3790 [math.CA]} .

\bibitem[{\citenamefont{{National Imagery and Mapping
  Agency}}(2000)}]{NIMA8350}
\bibinfo{author}{\bibnamefont{{National Imagery and Mapping Agency}}},
  \bibinfo{year}{2000}, \emph{\bibinfo{title}{Department {O}f {D}efense {W}orld
  {G}eodetic {S}ystem 1984}}, \bibinfo{type}{Technical Report}
  \bibinfo{number}{TR8350.2}, \bibinfo{institution}{{NIMA}},
  \urlprefix\url{http://earth-info.nga.mil/GandG/publications/tr8350.2/tr8350_2.html}.

\bibitem[{\citenamefont{\"Olander}(1952)}]{OlanderBullGeo26}
\bibinfo{author}{\bibnamefont{\"Olander}, \bibfnamefont{V.~R.}},
  \bibinfo{year}{1952}, \bibinfo{journal}{Bull. Geod.}
  \textbf{\bibinfo{volume}{26}}(\bibinfo{number}{3}), \bibinfo{pages}{337}.

\bibitem[{\citenamefont{Pollard}(2002)}]{PollardJG76}
\bibinfo{author}{\bibnamefont{Pollard}, \bibfnamefont{J.}},
  \bibinfo{year}{2002}, \bibinfo{journal}{J. Geod.}
  \textbf{\bibinfo{volume}{76}}(\bibinfo{number}{1}), \bibinfo{pages}{36}.

\bibitem[{\citenamefont{Rapp}(1991)}]{RappGG1}
\bibinfo{author}{\bibnamefont{Rapp}, \bibfnamefont{R.~H.}},
  \bibinfo{year}{1991}, \emph{\bibinfo{title}{Geometric Geodesy Part 1}}
  (\bibinfo{publisher}{Ohio State University}, \bibinfo{address}{Columbus,
  Ohio}).

\bibitem[{\citenamefont{Reichardt}(1957)}]{Reichardt}
\bibinfo{author}{\bibnamefont{Reichardt}, \bibfnamefont{H.}},
  \bibinfo{year}{1957}, \emph{\bibinfo{title}{Vorlesungen \"uber Vektor- und
  Tensorrechnung}}, volume~\bibinfo{volume}{34} of
  \emph{\bibinfo{series}{Hochschulb\"ucher f\"ur Mathematik}}
  (\bibinfo{publisher}{Deutscher Verlag der Wissenschaften},
  \bibinfo{address}{Berlin}).

\bibitem[{\citenamefont{Saito}(1970)}]{SaitoBullGeo44}
\bibinfo{author}{\bibnamefont{Saito}, \bibfnamefont{T.}}, \bibinfo{year}{1970},
  \bibinfo{journal}{Bull. Geod.}
  \textbf{\bibinfo{volume}{44}}(\bibinfo{number}{4}), \bibinfo{pages}{341}.

\bibitem[{\citenamefont{Sj\"oberg}(2012)}]{SjobergJGS2}
\bibinfo{author}{\bibnamefont{Sj\"oberg}, \bibfnamefont{L.~E.}},
  \bibinfo{year}{2012}, \bibinfo{journal}{J. Geod. Sci.}
  \textbf{\bibinfo{volume}{2}}(\bibinfo{number}{3}), \bibinfo{pages}{162}.

\bibitem[{\citenamefont{Smart}(1949)}]{Smart}
\bibinfo{author}{\bibnamefont{Smart}, \bibfnamefont{W.~M.}},
  \bibinfo{year}{1949}, \emph{\bibinfo{title}{Text-book on Spherical
  Astronomy}} (\bibinfo{publisher}{Cambridge University Press},
  \bibinfo{address}{Cambridge}), \bibinfo{edition}{4th} edition.

\bibitem[{\citenamefont{Sodano}(1958)}]{SodanoJGeod32}
\bibinfo{author}{\bibnamefont{Sodano}, \bibfnamefont{E.~M.}},
  \bibinfo{year}{1958}, \bibinfo{journal}{Bull. Geod.}
  \textbf{\bibinfo{volume}{32}}(\bibinfo{number}{2}), \bibinfo{pages}{13}.

\bibitem[{\citenamefont{Thomas}(1965)}]{ThomasJGR70}
\bibinfo{author}{\bibnamefont{Thomas}, \bibfnamefont{P.~D.}},
  \bibinfo{year}{1965}, \bibinfo{journal}{J. Geophys.\ Res.}
  \textbf{\bibinfo{volume}{70}}(\bibinfo{number}{14}), \bibinfo{pages}{3331}.

\bibitem[{\citenamefont{Tienstra}(1951)}]{TienstraBullGeo25}
\bibinfo{author}{\bibnamefont{Tienstra}, \bibfnamefont{J.~M.}},
  \bibinfo{year}{1951}, \bibinfo{journal}{Bull. Geod.}
  \textbf{\bibinfo{volume}{25}}(\bibinfo{number}{1}), \bibinfo{pages}{7}.

\bibitem[{\citenamefont{Tseng}(2014)}]{TsengJN67}
\bibinfo{author}{\bibnamefont{Tseng}, \bibfnamefont{W.-K.}},
  \bibinfo{year}{2014}, \bibinfo{journal}{J. Navig.}
  \textbf{\bibinfo{volume}{67}}(\bibinfo{number}{5}), \bibinfo{pages}{825}.

\bibitem[{\citenamefont{Vermeille}(2004)}]{VermeilleJG78}
\bibinfo{author}{\bibnamefont{Vermeille}, \bibfnamefont{H.}},
  \bibinfo{year}{2004}, \bibinfo{journal}{J. Geod.}
  \textbf{\bibinfo{volume}{78}}(\bibinfo{number}{1--2}), \bibinfo{pages}{94}.

\bibitem[{\citenamefont{You}(1999)}]{YouQV}
\bibinfo{author}{\bibnamefont{You}, \bibfnamefont{R.-J.}},
  \bibinfo{year}{1999}, in \emph{\bibinfo{booktitle}{Quo vadis
  geodesia\ldots?}}, edited by
  \bibinfo{editor}{\bibfnamefont{F.}~\bibnamefont{Krumm}} and
  \bibinfo{editor}{\bibfnamefont{V.~S.} \bibnamefont{Schwarze}}
  (\bibinfo{publisher}{Universit\"at Stuttgart}), Technical Reports of Geodesy
  and GeoInformatics, ISSN \bibinfo{issn}{0933-2839}.

\bibitem[{\citenamefont{Zhang} \emph{et~al.}(2005)\citenamefont{Zhang, Hsu, Wu,
  Li, Wang, Chai, and Du}}]{ZhangJG79}
\bibinfo{author}{\bibnamefont{Zhang}, \bibfnamefont{C.-D.}},
  \bibinfo{author}{\bibfnamefont{H.~T.} \bibnamefont{Hsu}},
  \bibinfo{author}{\bibfnamefont{X.~P.} \bibnamefont{Wu}},
  \bibinfo{author}{\bibfnamefont{S.~S.} \bibnamefont{Li}},
  \bibinfo{author}{\bibfnamefont{Q.~B.} \bibnamefont{Wang}},
  \bibinfo{author}{\bibfnamefont{H.~Z.} \bibnamefont{Chai}}, and
  \bibinfo{author}{\bibfnamefont{L.}~\bibnamefont{Du}}, \bibinfo{year}{2005},
  \bibinfo{journal}{J.\ Geod.}
  \textbf{\bibinfo{volume}{79}}(\bibinfo{number}{8}), \bibinfo{pages}{413}.

\end{thebibliography}

\end{document}